\documentclass[11pt,reqno]{amsart}
\usepackage{euscript,amssymb}
\usepackage[arrow,matrix,curve]{xy}
\usepackage{array,longtable}
\usepackage{a4wide}

\newenvironment{m-theorem}{%
\vskip5pt\refstepcounter{stff}\trivlist \itemindent 0pt
\item[\hskip\labelsep\bf Theorem \thestff]%
\it\ignorespaces}{\endtrivlist\vskip5pt}%
\newenvironment{m-proposition}{%
\vskip5pt\refstepcounter{stff}\trivlist \itemindent 0pt
\item[\hskip\labelsep\bf Proposition \thestff]%
\it\ignorespaces}{\endtrivlist\vskip5pt}%
\newenvironment{m-corollary}{%
\vskip5pt\refstepcounter{stff}\trivlist \itemindent 0pt
\item[\hskip\labelsep\bf Corollary \thestff]%
\it\ignorespaces}{\endtrivlist\vskip5pt}%
\newenvironment{m-lemma}{%
\vskip5pt\refstepcounter{stff}\trivlist \itemindent 0pt
\item[\hskip\labelsep\bf Lemma \thestff]%
\it\ignorespaces}{\endtrivlist\vskip5pt}%
\newenvironment{m-definition}{%
\vskip5pt\refstepcounter{stff}\trivlist \itemindent 0pt
\item[\hskip\labelsep\bf Definition \thestff]%
\ignorespaces}{\endtrivlist\vskip5pt}%
\newenvironment{m-notation}{%
\vskip5pt\refstepcounter{stff}\trivlist \itemindent 0pt
\item[\hskip\labelsep\bf Notation \thestff]%
\ignorespaces}{\endtrivlist\vskip5pt}%
\newenvironment{m-example}{%
\vskip5pt\refstepcounter{stff}\trivlist \itemindent 0pt
\item[\hskip\labelsep\bf Example \thestff]%
\ignorespaces}{\endtrivlist\vskip5pt}
\newenvironment{m-remark}{%
\vskip5pt\refstepcounter{stff}\trivlist \itemindent 0pt
\item[\hskip\labelsep\bf Remark \thestff]%
\ignorespaces}{\endtrivlist\vskip5pt}
\newenvironment{m-question}{%
\vskip5pt\refstepcounter{stff}\trivlist \itemindent 0pt
\item[\hskip\labelsep\bf Question.]%
\ignorespaces}{\endtrivlist\vskip5pt}%
\newenvironment{thm-nono}{
\vskip5pt\trivlist \itemindent 0pt
\item[\hskip\labelsep\bf Theorem.]%
\it\ignorespaces}{\endtrivlist\vskip5pt}%
\newenvironment{m-thank}{%
\vskip5pt\trivlist \itemindent 0pt
\item[\hskip\labelsep\it Acknowledgments]%
\ignorespaces}{\endtrivlist\vskip5pt}%

\let\mt\mapsto

\font\tenmsa=msam10 %
\newcommand\hdashpiece{%
{\vrule height2.75pt depth-2.35pt width2.3pt \kern1.7pt}}%
\newcommand\hdashpieces{%
{\hdashpiece\hdashpiece\hdashpiece\hdashpiece}}%
\let\dashto\dashrightarrow
\newcommand\dashar{\mathrel{%
\hdashpieces\kern-0.4pt\hbox{\tenmsa K}}}%

\let\euf\EuScript 
\let\cal\mathcal
\let\mbb\mathbb

\DeclareFontFamily{OT1}{rsfs}{}
\DeclareFontShape{OT1}{rsfs}{n}{it}{<->rsfs10}{}
\DeclareMathAlphabet{\crl}{OT1}{rsfs}{n}{it}

\let\ovl\overline
\let\unl\underline
\let\unbar\underbar

\let\tld\tilde

\let\wht\widehat
\let\nit\noindent

\let\disp\displaystyle
\let\srel\stackrel

\let\lan\langle
\let\ran\rangle
\newcommand\lran[1]{{\lan #1\ran}}

\newcommand\Aut{\operatorname{\textrm{Aut}\kern1pt}}
\newcommand\cAut{\operatorname{\mathcal{A}\kern-1pt\textit{ut}\kern1pt}}
\newcommand\End{\operatorname{\rm{End}\kern1pt}}
\newcommand\cEnd{\operatorname{\mathcal{E}\kern-1pt\textit{nd}\kern1pt}}

\newcommand\cHom{\operatorname{\mathcal{H}\kern-1pt\textit{om}\kern1pt}}
\newcommand\Img{{\rm Im}}

\newcommand\Lie{\mathop{\rm Lie}\nolimits}

\newcommand{\Weyl}{\mathop{\rm{Weyl}}\nolimits}

\newcommand\Spec{\mathop{\rm Spec}\nolimits}

\newcommand\invq{{\slash\kern-2.5pt\slash}}

\numberwithin{equation}{section}
\numberwithin{figure}{section} 

\let\l\lambda

\let\si\sigma
\let\Si\Sigma

\let\sm\setminus

\newcommand\bbk{\mbox{\rm I\kern-1.5pt k}}
\newcommand\sbbk{\hbox{\scriptsize I{\kern-.8pt}k}}

\newcommand\bbR{{\mbb R}}

\newcommand\bone{{1\kern-0.57ex\rm l}}

\newcommand\eI{{\euf I}}

\newcommand\eN{{\euf N}}

\newcommand\eO{{\euf O}}

\newcommand\codim{{\rm codim}}

\renewcommand\det{{\rm det}}

\let\ges\geqslant
\let\les\leqslant

\newcommand{\cd}{\mathop{\rm cd}\nolimits}
\newcommand{\cExt}{\mathop{{\crl E}{\kern-1pt}xt}\nolimits}
\newcommand{\BB}{{BB}}
\newcommand{\Bht}{{\rm Bruhat}}
\newcommand{\kk}{{\Bbbk}}

\title[Cohomological dimension of stratified varieties] %
{About the cohomological dimension\\ of certain stratified varieties}
\author{Mihai Halic, Roshan Tajarod}
\address{}

\begin{document}

\begin{abstract}
We determine an upper bound for the cohomological dimension of the complement of 
a closed subset in a projective variety which possesses an appropriate stratification. 
We apply the result to several particular cases, including the Bialynicki-Birula 
stratification; in this latter case, the bound is optimal.
\end{abstract}

\maketitle


\section*{Introduction}

The cohomological dimension measures at what extent a variety is affine or complete. 
Affine varieties are characterized by the fact that their cohomological dimension vanishes 
(Serre's criterion). 
For an irreducible quasi-projective variety $W$, $\cd(W)\les\dim(W)$ (cf. \cite{groth}), 
and the inequality is strict unless $W$ is proper 
(cf. Lichtenbaum's criterion of properness \cite{hart-cdav,lyub}). 
So, for the complement $W$ of a non-empty closed subset of a projective variety $X$, 
holds $\cd(W)<\dim(W)$. 

In general, it is not easy to determine the cohomological dimension of non-complete varieties. 
This might be the reason why there are not so many articles in the literature which contain explicit computations: see \cite{hart-cdav,ogus,lyub,walt}. 
A quasi-projective variety $W$ is by definition the complement of a closed subset $Z$ in 
a projective variety $X$, so it is reasonable to expect that the properties of $Z$ have influence 
on $\cd(W)$ (e.g. if $Z$ is an ample divisor, then $\cd(W)=0$). 
This issue was investigated in \cite{hart-cdav} for the complement of subvarieties with ample 
normal bundle and in \cite{ottm} for ample subvarieties. 

The goal of this article is to compute the cohomological dimension of non-complete varieties 
which admit appropriate stratifications; the definition is inspired from \cite{ca-go}. 
The main result is Theorem \ref{thm:cd-lci}, which is an upper bound for the cohomological 
dimension of the varieties admitting affine bundle stratifications. We also discuss when is reached 
this upper bound. 

The best known examples arise in the context of multiplicative group actions on projective varieties, 
when the stratification we are referring to is known as the Bialynicki-Birula decomposition; in this 
case the upper bound of the theorem is optimal. This is discussed in the section \ref{ssct:sink}. 
Zero loci of sections in globally generated vector bundles is another class of examples, 
not related to group actions, which yield affine bundle stratifications; this is discussed in 
the section \ref{ssct:0loci}. Throughout this note, the varieties are defined over an algebraically 
closed field $\kk$. 


\section{Background material}\label{sct:setup} 

\begin{m-definition}\label{def:cd}(cf. \cite[pp. 170]{groth}) 
Let $W$ be a quasi-projective scheme. The \emph{cohomological dimension} of $W$, denoted by 
$\cd(W)$,  is the least integer $n$ such that, for any quasi-coherent sheaf $\crl F$ on $W$, holds: 
$$
H^t(W,\crl F)=0,\;\forall\,t>n.
$$
\end{m-definition}
Although we are interested in the cohomology of \emph{coherent} sheaves, the necessity to apply  
Leray-type spectral sequences for affine morphisms (see below) forces us to work in the quasi-coherent 
setting. 

\begin{m-lemma}\label{lm:affine}
If $f:W\to B$ is an affine morphism, then $\cd(W)\les\cd(B)$.
\end{m-lemma}
There are examples where the inequality is strict (cf. \ref{expl:<}). 
\begin{proof}
This is a direct consequence of \cite[Ch. III, Corollaire 1.3.3]{groth-ega3}: since $f$ is affine, for any quasi-coherent sheaf $\crl F$ on $W$, holds: $H^t(W,\crl F)\cong H^t(B,f_*\crl F)$, for all $t\ges0$. 
\end{proof}

\begin{m-definition}\label{def:decomp}
Let $X$ be a quasi-projective variety. The data $(Z_\bullet,Y_\bullet)$ 
consisting of a filtration by \emph{closed} subsets 
$$
\emptyset=Z_0\subset Z_1\subset\ldots\subset Z_r=X, 
$$ 
and a collection $\{Y_1,\ldots,Y_r\}$ of quasi-projective varieties is called 
a \emph{lci affine bundle stratification} if, for all $j=1,\ldots,r$, the following properties hold: 
\begin{enumerate}
\item $Z_j\sm Z_{j-1}\subset X$ is a locally complete intersection (lci for short), for $j\ges2$;
\item 
 there is an affine morphism  $f_j:Z_j\sm Z_{j-1}\to Y_j$.
\end{enumerate}
\end{m-definition}
The role of the varieties $Y_\bullet$ may not be clear at this point. \textit{A priori}, one could 
take $Y_j=Z_j\sm Z_{j-1}$ and $f_j$ the identity; however, the interest is to have $Y_j$ as low 
dimensional as possible. The definition is inspired from \cite{ca-go}, but it is not the same.


\section{The main result}

Throughout this section we assume that $(Z_\bullet,Y_\bullet)$ is a lci affine bundle stratification 
of a quasi-projective variety $X$. The local cohomology groups and sheaves are usually defined for 
\emph{closed} supports; here we work with \emph{locally closed} supports, 
defined in \cite[Expos\'e I]{groth-sga2}; for an algebraic introduction, the reader may consult 
\cite{brod-sharp}. 
The letters `$H$' and `$\crl H$' will stand for the local cohomology groups and sheaves respectively.  
We start by recalling the following well-known result 
(mentioned en passant in \cite[pp. 30, bottom]{groth-sga2}): 

\begin{m-lemma}\label{lm:lcd}
Let $A$ be an affine scheme and $Y\subset A$ be a $\delta$-codimensional, 
complete intersection subscheme. 
Then, for any quasi-coherent sheaf $\crl F$ on $A$, holds $\,H^t_Y(A,\crl F)=0$, for all $t>\delta$. 
Consequently, if $W$ is locally closed, lci in a quasi-projective scheme $X$, then 
$\crl H^t_W(X,\crl F)=0$, for all quasi-coherent sheaves $\crl F$ on $X$ and $t>\codim_X(W)$. 
\end{m-lemma}

\begin{proof}
It is done by induction on $\delta$; for $\delta=0$, the statement is Serre's criterion. 
For the inductive step, write $Y=Y_1\cap\ldots\cap Y_\delta$, where all $Y_j\subset A$ 
are hypersurfaces, and consider the exact sequence: 
$$
\ldots\to 
H^{t-1}_{Y_1\cap\ldots\cap Y_{\delta-1}\cap (A\sm Y_\delta)}(A\sm Y_\delta,\crl F)
\to
H^t_{Y_1\cap\ldots\cap Y_{\delta}}(A,\crl F)
\to
H^t_{Y_1\cap\ldots\cap Y_{\delta-1}}(A,\crl F)
\to\ldots
$$
Since both $A$ and $A\sm Y_\delta$ are affine, the induction hypothesis implies that 
the left- and the right-hand-side vanish, for $t-1>\delta-1$.  
\end{proof}

\begin{m-lemma}\label{lm:loccoh}
Let the situation as in \ref{def:decomp}. 
For any coherent sheaf $\crl F$ on $X$ holds: 
$$
H^t_{Z_j\sm Z_{j-1}}(X;\crl F)=0,\;\forall\,t>\cd(Y_j)+\codim_X(Z_j\sm Z_{j-1}),
\;\forall\,j=2,\ldots,r.
$$
\end{m-lemma}

\begin{proof}
To simplify the notation, we denote $W:=Z_j\sm Z_{j-1}$, $U:=X\sm Z_{j-1}$, $B:=Y_j$. 
Since $W\subset U$ is closed, we deduce that $H^t_{W}(X,\crl F)=H^t_{W}(U,\crl F)$. 
The local cohomology groups can be computed by means of a spectral sequence 
(cf. \cite[Expos\'e I, Th\'eor\`eme 2.6]{groth-sga2}):
$$
H^b(U,\crl H^a_W(\crl F))\;\Rightarrow\;H^{b+a}_W(U,\crl F).
$$
First, since $W\subset U$ is lci, the lemma \ref{lm:lcd} implies that $\crl H^a_W(\crl F)=0$, 
for all $a>\codim(W)$. Second, 
$\crl H^a_{W}(\crl F)=\underset{m}{\varinjlim}\,\cExt^a(\eO_{U}/\eI_{W}^m,\crl F)$ 
(cf. \cite[Expos\'e II, Th\'eor\`eme 2]{groth}), which implies that 
$$
H^b(U,\crl H^a_W(\crl F))=
\underset{m}{\varinjlim}\,H^b(U,\cExt^a(\eO_{U}/\eI_{W}^m,\crl F)).
$$
The $\cExt$ groups are supported on (thickenings of) $W$ and $\cd(W)\les\cd(Y)$, 
so the expression above vanishes for $b>\cd(Y)$. 
The spectral sequence yields $H^t_W(U,\crl F)=0$ for $t$ in the indicated range.  
\end{proof}

\begin{m-theorem}\label{thm:cd-lci}
Let $(Z_\bullet,Y_\bullet)$ be a lci affine bundle stratification of a quasi-projective variety $X$. 
Then holds: 
\begin{equation}\label{eq:ineq-cd}
\cd(X\sm Z_1)\les\max_{j=2,\ldots,r}\{\cd(Y_j)+\codim_X(Z_j\sm Z_{j-1})\}.
\end{equation}
\end{m-theorem}
As we will see, in many geometric situations, the inequality above is actually an equality. 
However, the inequality is strict in general (cf. example \ref{expl:<}). 
\begin{proof}
Let $\crl F$ be an arbitrary quasi-coherent sheaf on $X\sm Z_1$; denote by $t_0$ the right hand side 
above. We show by decreasing recurrence that $H^t(X\sm Z_j,\crl F)=0$, for $j=r-1,\ldots,1$ and 
$t>t_0$. 
For $j=r-1$, we have 
$H^t(X\sm Z_{r-1},\crl F)\srel{\ref{lm:affine}}{=}H^t(Y_r,f_{r*}\crl F)=0$, for $t>\cd(Y_r)$. 
Now we prove the recursive step: the exact sequence in local cohomology 
(cf. \cite[Expos\'e I, Th\'eor\`eme 2.8]{groth-sga2}) yields: 
$$
\ldots\to H^t_{Z_j\sm Z_{j-1}}(X\sm Z_{j-1},\crl F)
\to H^t(X\sm Z_{j-1},\crl F)\to H^t(X\sm Z_j,\crl F)\to\ldots
$$
By the lemma \ref{lm:loccoh}, the left hand side vanishes for $t>t_0$, while the right hand side 
vanishes by assumption; thus the cohomology vanishing holds for $j-1$ too. 
After $r-1$ steps, we obtain the conclusion. 
\end{proof}

Henceforth we assume that $X,Y_\bullet$ are projective. Some natural questions occur: 
\begin{center}
Is \eqref{eq:ineq-cd} an equality? Under which assumptions is the equality reached?
\end{center}
Assume that the right-hand-side of \eqref{eq:ineq-cd} attains its maximum at 
$j_0\in\{2,\ldots,n\}$. Clearly, the answer to the question is affirmative if there is a 
$\big(\dim Y_{j_0}+\codim(Z_{j_0}\sm Z_{j_0-1})\big)$-dimensional, 
irreducible \emph{projective} variety $Z'_{j_0}\subset X\sm Z_1$. Indeed, 
the inequality `$\ges$' holds too, by Lichtenbaum's criterion \cite[Corollary 3.2]{hart-cdav}: 
$\;\cd(X\sm Z_1)\ges\cd(Z'_{j_0})=\dim Z'_{j_0}. $

To understand how realistic this naive answer is, we analyse the right-hand-side of \eqref{eq:ineq-cd}. 
There is an open affine subset $A\subset X$ and an irreducible projective variety $\tld Y_j\subset Z_j$ (a multi-section of $f_j$) with the following properties: 

-- $\tld Y_j\cap (Z_j\sm Z_{j-1})\neq\emptyset$ is lci in $Z_j\sm Z_{j-1}$;

-- $f_j:\tld Y_j\dashto Y_j$ is generically finite, in particular $\dim\tld Y_j=\dim Y_j$, for all $j$. 

\begin{m-remark}\label{rmk:lift}
Suppose that $X$ is smooth. Then the following hold: 
\begin{enumerate}
\item[--] There is an open affine subset $A\subset X$ such that 
both $(Z_j\sm Z_j-1)\cap A$ and $\tld Y_j\cap A$ are smooth; 

\item[--] The right hand side of \eqref{eq:ineq-cd} equals 
$d_j:=\dim\big({\unl{\sf N}_{(Z_j\sm Z_{j-1})\cap A/A}\big|}_{\tld Y_j\cap A}\big)$, 
that is the restriction to $\tld Y_j\cap A$ of the \emph{total space} of the normal bundle 
of $Z_j\sm Z_{j-1}$ in $X$; 

\item[--] There is a closed subvariety $W_{j}\subset A$ which `extends' $\tld Y_j\cap A$ in the normal 
direction to $(Z_j\sm Z_{j-1})\cap A$, that is: 
$\;\dim W_{j}=d_j,\quad  W_{j}\cap (Z_j\sm Z_{j-1})\cap A=\tld Y_j\cap A.$ 
\end{enumerate}
\end{m-remark}

\begin{proof}
The first claim holds for $A$ sufficiently small. Then use  that 
$\unl{\sf N}_{(Z_j\sm Z_{j-1})\cap A/A}$ is locally free, of rank $\codim_X(Z_j\sm Z_{j-1})$. 

Finally, consider $A=\Spec(R)$ and $(Z_j\sm Z_{j-1})\cap A=\Spec(A/I)$, where $R$ is a regular local ring and $I\subset R$ is a complete intersection ideal. Then $I^n/I^{n+1}$ are free $R/I$-modules, for $n\ges1$, and one can successively lift (not canonically) the identity $R/I\to R/I$ to a a ring homomorphism $R/I\to\hat R_I$. This latter induces an isomorphism between the formal completions: 
$\;\psi:\frac{R}{I}\big[\!\big[\frac{I}{I^2}\big]\!\big]\srel{\cong}{\to}\wht{R}_I$. 

Let $\frac{J}{I}\subset\frac{R}{I}$ be the ideal defining the subvariety $\tld Y_j\cap A$; it defines 
the closed (for the $I$-adic topology) ideal 
$\psi\big(\frac{J}{I}\big[\!\big[\frac{I}{I^2}\big]\!\big]\big)\subset\hat A_I$. 
Thus it lifts to an ideal $J'\subset A$ with the desired properties.
\end{proof}

Consequently, if one can ensure that $Z'_{n-j+1}:=\ovl{W_j}^X$ is contained in $X\sm Z_1$ 
(intuitively, if $W_j$ does not `bend backwards'), then $Z'_j$ satisfies the conditions in the 
naive answer.  

\begin{m-corollary}\label{cor:=cd}
Let $X$ be a projective variety and $(Z_\bullet,Y_\bullet)$ be an lci affine bundle stratification of it, 
with $Y_\bullet$ projective. 
Assume that, for $j=2,\ldots,n$, there is $Z'_{n+1-j}\subset X\sm Z_{j-1}$ irreducible, closed, 
such that the following hold:
$$
Z'_{n+1-j}\cap Z_j=Y_j,\quad\dim Z'_{n+1-j}+\dim Z_j=\dim Y_j+\dim X.
$$
Then we have 
$\disp\cd(X\sm Z_1)=\max_{j=2,\ldots,r}\{\cd(Y_j)+\codim_X(Z_j\sm Z_{j-1})\}.$
\end{m-corollary}
\nit As we will see, the BB-stratification of a projective variety satisfies this condition.


\section{Applications}\label{sct:applic}



\subsection{The complement of the sink of a $G_m$-action}\label{ssct:sink}

The Bialynicki-Birula (BB for short) decomposition arises in the context of the actions 
of the multiplicative group $G_m$. 
Let $X$ be a smooth projective $G_m$-variety; we denote the action by $\l:G_m\times X\to X$, 
and assume that the action is effective. 
In this situation $X$ admits two (the plus and minus) decompositions into smooth, 
locally closed subsets (cf. \cite{bb}): 
\begin{itemize} 
	\item 
		For any $x\in X$, the specializations at $\{0,\infty\}=\mbb P^1\sm G_m$ are denoted 
		$\underset{t\to 0}{\lim}\l(t)\times x$ and $\underset{t\to\infty}{\lim}\l(t)\times x$; 
		they are both fixed by $\l$. 
	\item 
		The fixed locus $X^\l$ of the action is a disjoint union 
		$\underset{s\in S_\BB}{\coprod}\kern-1ex Y_s$ of smooth subvarieties. 
		For $s\in S_\BB$, 
		$Y_s^\pm:=\{x\in X\mid\underset{t\to 0\;\text{resp.}\,\infty}{\lim}\kern-2ex
		\l(t)\times x\in Y_s\}$ is locally closed in $X$ (a BB-cell).
	\item 
		$X=\underset{s\in S_\BB}{\coprod}\kern-1ex Y_s^+
		=\underset{s\in S_\BB}{\coprod}\kern-1ex Y_s^-$, and the morphisms 
		$Y_s^\pm\to Y_s,\; x\mt\underset{t\to0\;\text{resp.}\,\infty}{\lim}\kern-2ex\l(t)\times x$ 
		are locally trivial, affine space fibrations. They are not necessarily vector bundles, 
		that is the transition functions may be non-linear. 
	\item 
		The \emph{source} $Y_{\rm source}$ and the \emph{sink} $Y_{\text{\rm sink}}$ 
		of the action are characterized by the fact that 
		$Y^+_{\rm source}\subset X$ is open, $Y^-_{\rm source}=Y_{\rm source}$, and 
		$Y_{\rm sink}^+=Y_{\rm sink}$, $Y^-_{\rm sink}\subset X$ is open. 
	\item 
		There is a (not unique) partial order `$\prec$' on $S_\BB$ such that 
		$\ovl{Y^+_s}\subset\underset{t\preceq s}{\bigcup}Y^+_{t}=:Z_s$ (cf. \cite{bb-order}); 
		the difference of two consecutive terms of the filtration $Z_\bullet$ is of the form $Y^+_s$. 
		The minimal element of this (plus) filtration (the $Z_1$ in \ref{def:decomp}) is $Y_{\rm sink}$. 
		A similar statement holds for the minus-decomposition. 
\end{itemize}

In this case, the theorem \ref{thm:cd-lci} yields the following:

\begin{m-theorem}\label{thm:cd-bb}
Let the situation be as above. Then holds:
$$
\cd(X\sm Y_{\rm sink})=\max\{\dim Y^-_s\mid s\neq{\rm sink}\}
=\dim(X\sm Y^-_{\rm sink}),
$$
and similarly $\cd(X\sm Y_{\rm source})=\dim(X\sm Y^+_{\rm source})$. 
\end{m-theorem}

\begin{proof}
Indeed, for all $s\in S_\BB$, we have 
$\dim Y^+_s+\dim Y^-_s=\dim Y_s+\dim X$, and we apply the corollary \ref{cor:=cd}. 
\end{proof}

\nit In the paragraph \ref{sssct:homog} we will need a few more details about the BB-decomposition; 
define: 
\begin{equation}\label{eq:l}
\begin{array}{rl}
G(\l)&:=\big\{g\in G\mid g^{-1}\l(t)g=\l(t),\;\forall\,t\in G_m\big\}
\text{ the centralizer of $\l$ in $G$}, 
\\[1ex]
P(\pm\l)&\disp
:=\big\{g\in G\mid\lim_{t\to 0}
\big(\l(t)^{\pm}g\,\l(t)^{\mp}\big)\text{ exists in }G\;\big(\text{and belongs to }G(\l)\big)%
\big\},
\\[1ex] 
U(\pm\l)&\disp 
:=\big\{g\in G\mid\lim_{t\to 0}
\big(\l(t)^{\pm}g\,\l(t)^{\mp}\big)=e\in G\big\}.
\end{array}
\end{equation}
Then $G(\l)$ is a connected, reductive subgroup of $G$, $P(\pm\l)\subset G$ are 
parabolic subgroups,  $G(\l)$ is their Levi-component, and $U(\pm\l)$ the unipotent 
radical (cf. \cite[\S 13.4]{sp}). 

\begin{m-lemma}\label{lm:y-invar}
\nit{\rm(i)} $Y_{\rm source},Y_{\rm sink}$ are invariant under $P(-\l)$ and $P(\l)$ respectively. 

\nit{\rm(ii)} $Y_s^+$ is $P(\l)$-invariant and $U(\l)$ preserves the fibration 
$Y_s^+\to Y_s$, for all $s\in S_\BB$. 
\end{m-lemma}

\begin{proof} 
(i) We prove the statement for $Y_{\rm source}$; the case $Y_{\rm sink}$ is analogous. 
We claim that $G(\l)$ leaves $Y_{\rm source}$ invariant; 
for $y\in Y_{\rm source}$ and $c\in G(\l)$ holds: 
$$y=\big(c^{-1}\l(t)c\big)y\;\Rightarrow\;cy=\l(t)\cdot(cy)\;\forall\,t\in G_m
\;\Rightarrow\;cy\in X^\l;\quad\text{thus }G(\l) y\subset X^\l.$$
But $X^\l$ is the disjoint union of its components, $G(\l)y$ is connected, 
and contains $y\in Y_{\rm source}$, so $G(\l)Y_{\rm source}=Y_{\rm source}$. 
The same argument shows that $G(\l)Y_s=Y_s$, for any $s\in S_\BB$.
For $g\in P(-\l)$ holds $c:=\lim_{t\to0}\l(t)^{-1}g\l(t)\in G(\l)$, so: 
$$
cy=\lim_{t\to0}\l(t)^{-1}g\l(t)y\in Y_{\rm source}\;\Rightarrow\;
\lim_{t\to0}\l(t)^{-1}(gy)\in Y_{\rm source}. 
$$ 
We claim that $gy\in Y_{\rm source}$; otherwise $gy\in X\sm Y_{\rm source}$ 
is `repelled' from $Y_{\rm source}$ and the limit belongs to another component of $X^\l$. 

\nit(ii) Take $x\,{\in}\,X$ with 
$\underset{t\to 0}{\lim}\l(t)x\,{=}\,y\in Y_s\,{\subset}\,X^\l$, $g\in P(\l)$: then 
$\disp c\,{:=}\lim_{t\to 0}\l(t)g\l(t)^{-1}{\in}\,G(\l)$,  
$$
\lim_{t\to 0}\l(t)gx
{=}\,
\lim_{t\to 0}\!\big(\l(t)g\l(t)^{-1}\!\cdot\!\l(t)x\big){=}\,cy 
\;\Rightarrow\;
gx\in Y_s^+.
$$
Similarly, for $g\in U(\l)$, one may check that $\underset{t\to 0}{\lim}\l(t)x=y$ 
implies $\underset{t\to 0}{\lim}\l(t)(gx)=y$. 
\end{proof}


\subsubsection{The case of toric varieties}\label{sssct:toric}

Let $N$ be a lattice of rank $d$, $\Si\subset N_\bbR$ be a projective simplicial fan, and $X_\Si$ 
be the corresponding toric variety; denote by $\{D_\rho\}_{\rho\in\Si(1)}$ the invariant divisors. 
We denote by $T\subset X_\Si$ the big torus, and consider a 1-PS $\l$ of $T$; note that 
$\l$ is determined by an element of $N$, which will be denoted the same (so $\l\in N$). 
Further details about toric varieties can be found in \cite{oda}. 

The fixed components of the $G_m$-action on $X$ determined by $\l$ is a disjoint union of toric 
subvarieties of $X_\Si$ (they are all $T$-invariant), which are intersections of $T$-invariant divisors. 
Hence they are of the form 
$$Y_\si:=\kern-1ex\underset{\rho\in\si(1)}{\mbox{$\bigcap$}}\kern-1ex D_\rho,\quad\si\in\Si.$$ 
Since $\Si$ is simplicial, $Y_\si$ is a complete intersection. 

An open affine, $T$-invariant neighbourhood of the generic point of $Y_\si$ in $X_\Si$ is $\Spec(\kk[\si^\vee])$, where 
$\si^\vee:=\{m\in N^\vee\mid m(\xi)\ges0,\;\forall\,\xi\in\si(1)\}.$ Consider the cones: 
$$
\si^\perp:=\{m\in\si^\vee\mid m(\xi)=0,\;\forall\,\xi\in\si\},\quad 
\si^{>0}:=\{m\in\si^\vee\mid \exists\,\xi\in\si,\;m(\xi)>0\}.
$$
Then $\si^{>0}$ determines the ideal $\cal I\subset\kk[\si^\vee]$ with quotient ring $\kk[\si^\perp]$, which induces the inclusion $\Spec(\kk[\si^\perp])\subset\Spec(\kk[\si^\vee])$; the left hand side is an affine chart for $Y_\si$ (more precisely, its big torus).

\begin{m-lemma}\label{lm:toric}
The fixed components $X_\Si$ are 
$Y_\si=\kern-1ex\underset{\rho\in\si(1)}{\bigcap}\kern-1ex D_\rho$, 
where $\si\in\Si$ are \emph{the minimal cones} such that $\l\in\lran{\si(1)}$ 
(the vector space generated by $\si(1)$). 

Let $\si_{\rm source}$ (respectively $\si_{\rm sink}$) be the minimal cones such that 
$\l\in{\rm interior}(\si_{\rm source})$ (respectively $-\l\in{\rm interior}(\si_{\rm sink})$). 
Then we have 
$Y_{\rm source}=Y_{\si_{\rm source}}$ and $Y_{\rm sink}=Y_{\si_{\rm sink}}$. 
\end{m-lemma}

\begin{proof}
Let $\si\in\Si$ be such that $Y_\si$ is fixed by $\l$. 
The closed points of $\Spec(\kk[\si^\vee])$ correspond to semi-group homomorphisms 
$x:(\si^\vee,+)\to(\kk,\cdot);$ the distinguished element 
$$x_\si:\si^\vee\to\kk,\quad x_\si(m):=\bigg\{
\begin{array}{cl}
1,&\text{if }m\in\si^\perp;\\ 0,&\text{if }m\in\si^{>0}
\end{array}$$
corresponds to the generic point of $\Spec(\kk[\si^\perp])$. The action of $\l$ on $x_\si$ is the following: 
$$\big(\l(t)\cdot x_\si\big)(m)=t^{m(\l)}\cdot x_\si(m),\;\forall\,m\in\si^\vee.$$
By assumption, it holds $\l(t)\cdot x_\si=x_\si,\;\forall\,t\in G_m$. 
For $m\in\si^{>0}$, both sides vanish. For $m\in\si^\perp$, we deduce that $m(\l)=0$; since $m$ is arbitrary, it follows that $\l\in\lran{\si(1)}$. 
\end{proof}

Let $\Si_\BB\subset\Si$ be the subset consisting of the cones $\si$ as above; 
for $\rho\in\Si(1)$, denote by $\xi_\rho$ the generator of $\rho\cap N$. 
Then, for any $\si\in\Si_\BB$, we can (uniquely) write 
$$
\l=
\underset{\rho\in{\si(1)}^-}{\sum}\underbrace{c_\rho\,\cdot}_{<0}\xi_\rho
+
\underset{\rho\in{\si(1)}^+}{\sum}\underbrace{c_\rho\,\cdot}_{>0}\xi_\rho.
$$

\begin{m-corollary}\label{cor:toric}
Let the situation be as above. Then holds: 
$$\cd(X\sm Y_{\rm sink})=d-\min\{\;\#{\,\si(1)}^+\mid\si\in\Si_\BB\}.$$
\end{m-corollary}

\begin{proof}
Note that $\codim Y^-_s=\#{\,\si(1)}^+$, for all $\si\in\Si_\BB$, and apply 
the theorem \ref{thm:cd-bb}.
\end{proof}


\subsubsection{The case of homogeneous varieties}\label{sssct:homog}

Consider the homogeneous variety $X=G/P$, where $G$ is connected, reductive, and $P$ is 
a parabolic subgroup, and consider a subgroup $\l:G_m\to G$; it induces a $G_m$-action on $X$. 
The adjoint action of $\l$ on $\Lie(G)$ decomposes it into the positive/negative weight spaces:
$\;\Lie(G)=\Lie(G)^-_\l\oplus\Lie(G)^0\oplus\Lie(G)^+_\l.$

\begin{m-lemma}\label{lm:y-homog} 
The following statements hold:
\\ \nit{\rm(i)} 
The connected components of the fixed locus $X^\l$ are homogeneous for the action of $G(\l)$.
\\ \nit{\rm(ii)} 
The sink $\,Y_{\rm sink}$ contains $\hat e\in G/P$ if and only if 
$\l\subset P$ and $\Lie(G)^+_\l\subset\Lie(P)$. 
\end{m-lemma}

\begin{proof}
(i) The differential of the multiplication  
$\Lie(G)\to T_yX$ is surjective at any $y\in X^\l$, and is $\l$-equivariant 
for the adjoint action on $\Lie(G)$. Both sides decompose into direct sums of weight 
spaces; in particular, $\Lie(G(\l))=\Lie(G)^0_\l\to (T_yX)^0$ is surjective too. 
Therefore all the $G(\l)$-orbits are open in $X^\l$, hence the components of $X^\l$ 
are homogeneous under $G(\l)$. 

\nit(ii) 
The point $\hat e$ belongs to $Y_{\rm sink}$ if and only if: 

-- $\hat e$ is fixed by $\l$, that is $\Img(\l)\subset P$; 

-- the weights of $\l$ on $T_{\hat e}X\cong\Lie(G)/\Lie(P)$ are negative. 
\end{proof}

If $Q,P\subset G$ are two parabolic subgroups, $G/P$ decomposes into the following 
finite disjoint union (the Bruhat decomposition, cf. \cite[\S 8]{sp}) of locally closed orbits 
under the action of $Q$:
$$
G/P=\underset{w\in S_\Bht}{\coprod}\kern-1ex QwP,
\;\text{with}\;
S_\Bht=\Weyl(Q)\big\backslash\Weyl(G)\,\big\slash\Weyl(P).
$$
Actually, $S_\Bht$ parameterizes the $\Weyl(Q)$-orbits in $(G/P)^T$. 
Each double coset in $S_\Bht$ contains a unique representative of minimal length; 
for each $w\in S_\Bht$ of minimal length, 
$$
\dim (QwP)=
{\rm length}(w)+\dim\big({\rm Levi}(Q)\big/{\rm Levi}(Q)\cap wPw^{-1}\big).
$$

\begin{m-proposition}\label{prop:BB-homog}
The BB-decomposition of $G/P$ for the action of $\l$ coincides with 
the Bruhat decomposition for the action of $P(\l)$. \\ 
If $\Lie(G)^+_\l\subset P$, then the sink of the action is $P(\l)P\cong G(\l)/G(\l)\cap P$. 
\end{m-proposition}
We remark that any standard parabolic $Q\subset G$ occurs as $P(\l)$, for some $\l$ 
such that $\Lie(G)^+_\l$ is contained in the Borel subgroup of $\Lie(G)$. 
\begin{proof}
Each Bruhat cell is the $P(\l)$-orbit of some $x\in (G/P)^T$; 
any such $x$ belongs to a component $Y_s\subset (G/P)^\l$; 
finally, the \BB-cell $Y_s^+$ is $P(\l)$-invariant (cf. lemma \ref{lm:y-invar}). 
Hence each Bruhat cell is contained in a unique \BB-cell. 
But the union of the former is $G/B$, and the latter cells are pairwise disjoint. 
It follows that each Bruhat cell equals some \BB-cell. 
\end{proof}

\begin{m-theorem}\label{thm:bruhat}
Let the situation be as above. Then holds 
$$
\cd\bigg(\frac{G}{P}\sm \frac{G(\l)}{G(\l)\cap P}\bigg)
=\max_{w\in S_\Bht\sm\{e\}}\dim\big(\,P(-\l)wP\,\big).
$$
\end{m-theorem}

\begin{proof}
It is a direct consequence of the theorem \ref{thm:cd-bb}: 
the minus BB-decomposition corresponds to the action of $-\l$, 
so we must replace $P(\l)\rightsquigarrow P(-\l)$. 
\end{proof}


\subsection{Complements of zero loci in globally generated vector bundles}\label{ssct:0loci} 
The previous examples might have given the impression that lci affine bundle stratifications 
belong to the realm of group actions. Here we show that this is not the case. 

\begin{m-proposition}
Let $Z$ be a closed subscheme of a projective variety $X$. We assume that there is a modification 
$\tld X$ of $X$ along $Z$ such that the following holds: there is a projective variety $X'$ and 
a morphism $f:\tld X\to X'$ such that the exceptional divisor $E_Z$ is $f$-relatively ample. 
Then $\cd(X\sm Z)\les\dim X'$. 
\end{m-proposition}

\begin{proof}
The assumption that $E_Z$ is $f$-relatively ample implies that $f:\tld X\sm E_Z\to X'$ 
is an affine morphism. Since $\tld X\sm E_Z\cong X\sm Z$, in this case we obtain 
a filtration with two strata: $Z_1:=E_Z$ and $Z_2:=\tld X$. 
\end{proof}

The previous situation arises as follows. 
Consider a globally generated vector bundle $\eN$ of rank $\nu$ on a projective variety $X$, 
such that $\det(\eN)$ is ample. Let $Z$ be the zero locus of an arbitrary (non-zero) section in $\eN$. 

\begin{equation}\label{eq:tld-x}
\xymatrix@R=2em@C=1em{
\tld X\ar@{^(->}[r]\ar[d]_-\pi\ar@<-5pt>[rrd]_-f
&
\mbb P(\eN)
{=}\,
\mbb P\Bigl(
\mbox{$\overset{\nu-1}{\bigwedge}$}\eN^\vee\otimes\det(\eN)
\Bigr)\ar@{^(->}[r]
&
X\times\mbb P\Bigl(
\mbox{$\overset{\nu-1}{\bigwedge}$}\Gamma(\eN)^\vee
\Bigr)\ar[d]
\\ 
X&&\mbb P:=\mbb P\Bigl(
\mbox{$\overset{\nu-1}{\bigwedge}$}\Gamma(\eN)^\vee
\Bigr),
}
\end{equation}
and holds 
\begin{equation}\label{eq:o1}
\eO_{\tld X}(E_Z)=\eO_{\mbb P(\eN)}(-1)\big|_{\tld X}=
\bigl(\det(\eN)\boxtimes\eO_{\mbb P}(-1)\bigr)\big|_{\tld X},
\end{equation}
so $E_Z$ is $f$-relatively ample.

\begin{m-corollary}\label{cor:0sct}
Let the situation be as above. Then holds $\cd(X\sm Z)\les\dim f(\tld X)$.
\end{m-corollary}

\begin{m-example}\label{expl:<}
In general the previous inequality is strict. 
Let $\eN:=\eO_{\mbb P^3}(1)\oplus\eO_{\mbb P^3}(2)$, and let $Z\subset\mbb P^3$ 
be the zero locus of a general section in $\eN$: it is the intersection of a plane $\{s_1=0\}$ 
with a quadric $\{s_2=0\}$. On one hand, $Z\subset\mbb P^3$ is an ample, $2$-codimensional 
subvariety, so $\cd(\mbb P^3\sm Z)=1$  (cf. \cite[Theorem 7.1, 5.4]{ottm}). 

On the other hand, the morphism $f$ in \eqref{eq:tld-x} is the following: 
$$
\begin{array}{l}
\mbb P^3\dashto
\mbb P\big(\eO(1)_{\mbb P^3}\oplus\eO_{\mbb P^3}(2)\big)\cong
\mbb P\big(\eO_{\mbb P^3}(-1)\oplus\eO_{\mbb P^3}\big)\to\mbb P^4,
\\[1ex]  
\unbar{\mbox{$x$}}\mt[x_0s_1(x):x_1s_1(x):x_2s(x):x_3s_1(x):s_2(x)]. 
\end{array}
$$
We claim that the image of this morphism is $2$-dimensional. 
Indeed, consider 
$$\unbar{\mbox{$y$}}=[y_0:\ldots:y_4]=f(\unbar{\mbox{$x$}})\in\Img(f),\quad y_0=1.
$$ 
One computes 
$x_0=1,x_1=y_1,x_2=y_2,x_3=y_3, s_1(1,y_1,y_2,y_3)=1,s_2(1,y_1,y_2,y_3)=y_4.$ 
Thus $\unbar{\mbox{$y$}}$ satisfies two independent equations in $\mbb P^4$. 
\end{m-example}


\end{document}